\documentclass[11pt]{amsart}
\usepackage{color}
\usepackage{epsfig}
\usepackage{verbatim}

\begin{document}
\title {Strong Splitter Theorem}
\maketitle 
\begin {center}
S. R. Kingan 
\footnote{The first author is partially supported by  PSC-CUNY grant number 64181-00 42} \\     
Department of Mathematics \\
Brooklyn College, City University of New York\\
 Brooklyn, NY 11210\\
skingan@brooklyn.cuny.edu\\  
\bigskip

Manoel Lemos 
\footnote{The second author is partially supported by CNPq under grant number 300242/2008-05.}\\
Departamento de Matematica \\
Universidade Federal de Pernambuco\\
Recife, Pernambuco, 50740-540, Brazil\\
manoel@dmat.ufpe.br\\  
\end {center}
\bigskip

\begin{abstract} The Splitter Theorem states that, if $N$ is a $3$-connected proper minor of a 3-connected matroid $M$ such that, if $N$ is a wheel or whirl then $M$ has no larger wheel or whirl, respectively, then there is a sequence $M_0, \dots , M_n$ of $3$-connected matroids with $M_0\cong N$, $M_n=M$ and for $i\in \{1, \dots , n\}$, $M_i$ is a single-element extension or coextension of $M_{i-1}$. Observe that there is no condition on how many extensions may occur before a coextension must occur. In this paper, we give a strengthening of the Splitter Theorem, as a result of which we can obtain, up to isomorphism, $M$ starting with $N$ and at each step doing a 3-connected single-element extension or coextension, such that at most two consecutive single-element extensions occur in the sequence (unless the rank of the matroids involved are $r(M)$). Moreover, if two consecutive single-element extensions by elements $\{e, f\}$ are followed by a coextension by element $g$, then $\{e, f, g\}$ form a triad in the resulting matroid. Using the Strong Splitter Theorem, we make progress toward the problem of determining the almost-regular matroids [6, 15.9.8]. {\it Find all $3$-connected non-regular matroids such that, for all $e$, either $M\backslash e$ or $M/e$ is regular.} In [4] we determined the binary almost-regular matroids with at least one regular element (an element such that both $M\backslash e$ and $M/e$ is regular) by characterizing the class of binary almost-regular matroids with no minor isomorphic to one particular matroid that we called $E_5$. As a consequence of the Strong Splitter Theorem we can determine the class of binary matroids with an $E_5$-minor, but no $E_4$-minor.  
\end{abstract}

\bigskip 

\section {\bf Introduction}

The matroid terminology follows Oxley [6]. Let $M$ be a matroid and $X$ be a subset of the ground set $E$. The {\it connectivity function}  $\lambda$ is defined as $\lambda (X) = r(X) + r(E-X) - r(M)$.  Observe that $\lambda (X) = \lambda (E-X)$.  For $k\ge 1$, a partition $(A, B)$ of $E$ is called a $k$-separation if $|A|\ge k$, $|B|\ge k$, and $\lambda (A) \le k-1$.  When $\lambda (A)=k-1$, we call $(A, B)$ an {\it exact k-separation}.  When $\lambda (A)=k-1$ and $|A|=k$ or $|B|=k$, we call $(A, B)$ a {\it minimal exact k-separation}. For $n\ge 2$, we say $M$ is {\it n-connected} if $M$ has no $k$-separation for $k\le n-1$.  A matroid is {\it internally $n$-connected} if it is $n$-connected and has no non-minimal exact $n$-separations. In particular, a simple matroid is 3-connected if $\lambda (A)\ge 2$ for all partitions $(A, B)$ with $|A|\ge 3$ and $|B|\ge 3$. A 3-connected matroid is {\it internally $4$-connected}  if $\lambda (A)\ge 3$ for all partitions $(A, B)$ with $|A|\ge 4$ and $|B|\ge 4$. To eliminate trivial cases, we shall also assume that a 3-connected matroid has at least four elements.  

If $M$ and $N$ are matroids on the sets $E$ and $E \cup e$ where $e   \not\in E$, then $M$ is a  {\it single-element extension} of $N$ if $M \backslash e = N$,  and $M$ is a {\it single-element coextension} of $N$ if $M^*$ is a single-element extension of $N^*$. If $N$ is a 3-connected matroid, then an extension $M$ of $N$ is  3-connected provided $e$ is not in a 1- or 2-element circuit of $N$ and $e$ is not a coloop of $N$. Likewise, $M$ is a 3-connected coextension of $N$ if $M^*$ is a 3-connected extension of $N^*$. 

In 1966 Tutte proved that for a $3$-connected matroid $M$ that is not a wheel or a whirl, there exists an element $e\in E(M)$, such that either $M\backslash e$ or $M/e$ is $3$-connected [11].
In other words, if $M$ is a $3$-connected matroid, then $M$ has a $3$-connected proper minor $M'$ such that $|E(M)-E(M')|=1$, unless $M$ is a wheel or whirl. In the case that $M$ is a wheel or a whirl, for every $e$ there is an $f$ such that $M\backslash e/f$ is 3-connected. So Tutte's theorem can be restated as follow: If $M$ is a $3$-connected matroid, then $M$ has a $3$-connected proper minor $M'$ such that $|E(M)-E(M')|\le 2$.

In 1972 Brylawski proved that if $M$ is a 2-connected matroid with a proper 2-connected minor $N$, then there exists  $e\in E(M)-E(N)$ such that either $M\backslash e$ or $M/e$ is 2-connected and has $N$ as a minor  [2]. So when the matroid is 2-connected we can maintain 2-connectivity in $M\backslash e$ or $M/e$, as well as the presence of a certain minor. In 1980 and 1981 Seymour and Tan independently proved that, if $N$ is a 3-connected proper minor of a 3-connected matroid $M$ such that if $N$ is a wheel or whirl then $M$ has no larger wheel or whirl, respectively, then there exists  $e\in E(M)-E(N)$ such that either $M\backslash e$ or $M/e$ is $3$-connected and has $N$ as a minor [7. 8]. This is known as the Splitter Theorem. 
In other words, $M$ has a $3$-connected minor $M'$ with 
$|E(M)-E(M')|=1$ and having an $N$-minor, unless $M$ is a wheel or whirl, in which case $M$ has a $3$-connected proper minor $M'$ with $|E(M)-E(M')|=2$ and having an $N$-minor. A formal statement of the Splitter Theorem appears below [6, 12.1.2]. 

\bigskip
\noindent{\bf Theorem 1.1.} {\it Suppose $N$ is a connected, simple, cosimple proper minor of a $3$-connected matroid $M$ such that, if  $N$ is a wheel or whirl then $M$ has no larger wheel or whirl-minor, respectively. Then $M$ has a connected, simple, cosimple minor $M'$ and an element $e$ such that  $M'\backslash e$ or $M'/ e$ is isomorphic to $N$.} $\qed$
\bigskip

The next two results are reformulations of the Splitter Theorem and appear in [6, 12.1.3] and [6, 12.2.1]. Let $\mathcal M$ be a class of matroids closed under minors and isomorphism.  A {\it splitter} $N$ for $\mathcal M$ is a 3-connected matroid in $\mathcal M$ such that no $3$-connected matroid in $\mathcal M$ has $N$ as a proper minor. 

\bigskip
\noindent{\bf Corollary 1.2.} {\it Suppose $N$ is a $3$-connected matroid in $\mathcal M$ such that, if  $N$ is a wheel or whirl then it is the largest wheel or whirl in $\mathcal M$. Suppose further that every $3$-connected single-element extension and coextension of $N$ does not belong to $\mathcal M$. Then $N$ is a splitter for $\mathcal M$.}    $\qed$
\bigskip

\noindent Checking if a matroid is a splitter is a potentially infinite task. The above reformulation of the Splitter Theorem turns this into a finite task, that can be easily checked. 

Next, suppose $\mathcal M$ is defined as having a specific $3$-connected matroid in it, for example, the class of $3$-connected binary matroids with an $F_7$ or $F_7^*$-minor, but without an $M(W_4)$-minor. The third reformulation of the Splitter Theorem asserts that the entire class can be built up by performing single-element extensions and coextensions starting with the specified matroid and checking for the specified excluded minor(s). 

\bigskip
\noindent{\bf Corollary 1.3.} {\it Suppose $N$ is a $3$-connected proper minor of a $3$-connected matroid $M$ such that, if  $N$ is a wheel or whirl then $M$ has no larger wheel or whirl-minor, respectively. Then,  there is a sequence $M_0, \dots , M_n$ of $3$-connected matroids with $M_0\cong N$, $M_n=M$ and for $i\in \{1, \dots , n\}$, $M_i$ is a single-element extension or coextension of $M_{i-1}$.} $\qed$
\bigskip

Subsequently, a couple of variations of the Splitter Theorem have been developed. For the variations of the Splitter Theorem, it suffices to state the results in terms of $N$ being 3-connected since the general case when $N$ is connected, simple, and cosimple is covered by the original Splitter Theorem.   Observe that the minor $M'$ is not required to have a single-element deletion or contraction equal to $N$, but only to have such a minor isomorphic to $N$. There is a counterexample to show the stronger statement of equality does not hold [6, 12.1 Ex. 7]. Truemper [6, 12.3.2]  strengthened the conclusion by proving that, if $N$ is a $3$-connected proper minor of a $3$-connected matroid $M$, then $M$ has a $3$-connected minor $M'$ and an element $e$ such that  $co(M'\backslash e)=N$ or $si(M'/ e)=N$ and $|E(M')-E(N)|\le 3$. Bixby and Coullard gave another similar variant [6, 12.3.6].

Coullard and Oxley [6, 12.3.1] showed that the restriction on excluding wheels and whirls can be weakened, so that instead of applying to all such matroids, it applies only to the smallest $3$-connected wheels and whirls. They proved that if $N$ is a $3$-connected proper minor of a $3$-connected matroid $M$ that is not a wheel or a whirl and if $N\cong W^2$, then $M$ has no $W^3$-minor and if $N\cong M(W_3)$, then $M$ has no $M(W_4)$-minor, then $M$ has a $3$-connected minor $M'$ and an element $e$ such that  $M'\backslash e$ or $M'/ e$ is isomorphic to $N$.

We prove a new variant of the Splitter Theorem, the usefulness of which becomes apparent in the third section, where we prove structural results by applying it.

\bigskip
\noindent{\bf Theorem 1.4.} {\it
Suppose $N$ is a $3$-connected proper minor of a $3$-connected matroid $M$ such that,  if  $N$ is a wheel or whirl then $M$ has no larger wheel or whirl-minor, respectively. Further, suppose $m=r(M)-r(N)$. Then there is a sequence of $3$-connected matroids $M_0,M_1,\dots,M_n$, for some integer $n\ge m$, such that
\begin{enumerate}
\item[(i)] $M_0\cong N$;
\item[(ii)] $M_n=M$;
\item[(iii)] for $k\in\{1,2,\dots,m\}$, $r(M_k)-r(M_{k-1})=1$ and $|E(M_k)-E(M_{k-1})|\le 3$; and
\item[(iv)] for $m<k\le n$, $r(M_k)=r(M)$ and $|E(M_k)-E(M_{k-1})|=1$.
\end{enumerate}
Moreover, when $|E(M_k)-E(M_{k-1})|=3$, for some $1\le k\le m$, $E(M_k)-E(M_{k-1})$ is a triad of $M_k$.}
\bigskip

Thus we can obtain, up to isomorphism,  $M$ starting with $N$ and at each step doing a 3-connected single-element extension or coextension, such that at most two consecutive single-element extensions occur in the sequence (unless the rank of the matroids involved are $r$). Moreover, if two consecutive single-element extensions by elements $\{e, f\}$ are followed by a coextension by element $g$, then $\{e, f, g\}$ form a triad in the resulting matroid. Finally, note that we can replace the restrictions on $M$ and $N$ by the weaker restrictions on $M$ and $N$ given by Coullard and Oxley.

\section {Proof of the Strong Splitter Theorem}

In this section we give the proof of Theorem 1.4. Let us begin by proving a  key lemma. 
\bigskip

\noindent {\bf Lemma 2.1.} {\it 
Suppose $N$ is a $3$-connected proper minor of a $3$-connected matroid $M$ and  $r(M)=r(N)+1$. Then, either
\begin{enumerate}
\item[(i)] There is an element $e\in E(M)-E(N)$ such that $M\backslash e$ is $3$-connected and $N$ is a minor of $M\backslash e$; or
\item[(ii)] $|E(M)-E(N)|\le 3$.
\end{enumerate}
Moreover, when $|E(M)-E(N)|=3$, $T^*=E(M)-E(N)$ is triad of $M$ and $M\backslash T^*=N$.}
\bigskip

\noindent {\bf Proof.}
There is a set $A$ of elements of $M$ and an element $b$ of $M$ such that $b\not\in A$ and $N=M\backslash A/b$. If $A=\emptyset$, then (ii) follows. Assume that $A\not=\emptyset$. Choose $e\in A$. If $M\backslash e$ is 3-connected, then (i) follows. Suppose that $M\backslash e$ is not 3-connected. Let $\{X,Y\}$ be a 2-separation for $M\backslash e$. As $\{X\cap E(N),Y\cap E(N)\}$ is not a 2-separation for $N$, it follows that $\min\{|X\cap E(N)|,|Y\cap E(N)|\}\le 1$, say $|Y\cap E(N)|\le 1$. We do not loose generality by assuming that $Y$ is closed in $M$. Using this 2-separation, we can decompose $M\backslash e$ as the 2-sum of matroids $M_X$ and $M_Y$ such that $E(M_X)=X\cup z$ and $E(M_Y)=Y\cup z$, for some $z\not\in E(M)$. Observe that $N$ is isomorphic to a minor of $M_X$ because $|E(N)-X|\le 1$. In particular, $r(M_X)\ge r(N)$. But
\begin{equation*}
r(N)+1=r(M)=r(M\backslash e)=r(X)+r(Y)-1=r(M_X)+r(M_Y)-1\ge r(N)+1
\end{equation*}
(As $M\backslash e$ is simple, it follows that $r(M_Y)=r(Y)\ge 2$.) We must have equality along this display. Therefore
\begin{equation}\label{11.11.11.a}
\hbox{$r(M_X)=r(N)$ and $r(M_Y)=2$.}
\end{equation}
In particular,
\begin{equation}\label{11.11.11.b}
b\in Y.
\end{equation}
We can be more precise  about a minor $N'$ of $M_X$ isomorphic to $N$. Observe that $N'=N=M_X\backslash [(A\cap X)\cup z]$, when $E(N)\cap Y=\emptyset$, or $N$ is obtained from $N'=M_X\backslash (A\cap X)$ by relabeling $z$ by $z'$, when $E(N)\cap Y=\{z'\}$. As $Y$ is closed in $M$, it follows that $M_X$ is simple because the series class of $z$ in $M_X$ is trivial. Observe that $M_X$ is 3-connected because $N'$ is a 3-connected restriction of $M_X$ having the same rank as $M_X$ and $M_X$ is simple. If $M_Y$ is not 3-connected, then $M_Y$ is not simple because $r(M_Y)=2$. In this case, there is an element $f\in Y$ such that $\{f,z\}$ is a parallel class of $M_Y$ and $M_Y\backslash f$ is 3-connected (recall that no two elements of $Y$ are in parallel in $M$ and so in $M_Y$). We have two possibilities to consider:
\begin{enumerate}
\item[(a)] $M_Y$ is 3-connected. In this case, $\{X,Y\}$ is the unique 2-separation for $M\backslash e$.
\item[(b)] $M_Y$ is not 3-connected. In this case, $\{X,Y\}$ and $\{X\cup f,Y-f\}$ are the  2-separations for $M\backslash e$.
\end{enumerate}

Assume that $g\in A\cap X$. As $N'$ and $M_X$ are 3-connected matroid with the same rank and $N'$ is a restriction of both $M_X$ and $M_X\backslash g$, it follows that $M_X\backslash g$ is 3-connected. Therefore $\{X-g,Y\}$ is the unique 2-separation for $M\backslash\{e,g\}$, when (a) occurs, or $\{X-g,Y\}$ and $\{(X-g)\cup f,Y-f\}$ are the 2-separations for $M\backslash\{e,g\}$, when (b) occurs. Moreover, $M\backslash\{e,g\}$ has no 1-separation. In both cases, each set in these 2-separations does not span $e$ and so $M\backslash g$ is 3-connected. We have (i) because $N$ is a minor of $M\backslash g$. We may assume that
\begin{equation}\label{09.11.11.a}
A\cap X=\emptyset.
\end{equation}
With a similar argument, we conclude that $M\backslash h$ is a 3-connected matroid having $N$ as a minor, when $h\in (Y-P)\cap A$, provided $|Y-P|\ge 3$, where $P$ is the parallel class of $M_Y$ containing $z$. (Since at most one element of $Y-P$ belongs to $E(N)$ and at most one element of $Y-P$ is equal to $b$, it follows that $Y-P$ meets $A$ provided $|Y-P|\ge 3$, that is, the element $h$ exists.) Thus (i) also occurs unless
\begin{equation}\label{09.11.11.b}
|Y-P|=2.
\end{equation}
We may assume the last identity otherwise the result follows. If $|P|=2$, that is, $P=\{z,f\}$, then $\{X,Y-f\}$ is the unique 2-separation for $M\backslash\{e,f\}$. Moreover, $f\in A$ or $f\in E(N)$. If $f\in A$, $M\backslash f$ is a 3-connected matroid having $N$ as a minor because $e$ is not spanned by any set in this 2-separation. In this case, we have (i). Hence, we could also assume that
\begin{equation}\label{09.11.11.c}
\hbox{if $|P|=2$, then $P-z\subseteq E(N)$.}
\end{equation}

By~\eqref{11.11.11.b},~\eqref{09.11.11.a} and~\eqref{09.11.11.c}, $E(M)-E(N)\subseteq (Y-P)\cup e$. By~\eqref{09.11.11.b}, $|E(M)-E(N)|\le 3$. The first part of this result follows. Now, we establish the second part. Assume that
\begin{equation*}
|E(M)-E(N)|=3.
\end{equation*}
In this case, no element of $Y-P$ belongs to $N$. So $T^*=e\cup (Y-P)=E(M)-{\rm cl}_M(X)$ is a triad of $M$ avoiding $E(N)$. But $b\in T^*$ and $T^*-b=A$. Therefore $N=M\backslash T^*$. $\qed$

\bigskip

\noindent {\bf Proof of Theorem 1.4. } Suppose $N$ is a $3$-connected proper minor of a $3$-connected matroid $M$ such that,  if  $N$ is a wheel or whirl then $M$ has no larger wheel or whirl-minor, respectively. By Corollary 1.3, there is a sequence of 3-connected matroids $N_0,N_1,\dots,N_n$ such that $N_0\cong N$, $N_n=M$ and, for each $i$ belonging to $\{1,2,\dots,n\}$, $N_i$ is a single-element extension or coextension of $N_{i-1}$. Set $m=r(M)-r(N)$. There are indexes $i_1,i_2,\dots,i_m$ such that $0<i_1<i_2<\dots<i_m$ and, for each $k$ belonging to $\{1,2,\dots,m\}$, 
\begin{equation*}
r(N_{i_k})-r(N_{i_k-1})=1
\end{equation*}
(That is, $N_{i_k}$ is the first matroid in the sequence $N_0,N_1,N_2,\dots,N_n$ having rank equal to $r(N)+k$.) Choose this sequence such that $(i_1,i_2,\dots,i_m)$ is minimal in the lexicographic order. By Lemma 2.1, $i_1=|E(N_{i_1})-E(N_0)|\le 3$ and, for each $k$ such that $2\le k\le m$, $i_k-i_{k-1}=|E(N_{i_k})-E(N_{i_{k-1}})|\le 3$. Moreover, when the equality holds, we have that $E(N_{i_1})-E(N_0)$ is a triad of $E(N_3)$, when $k=1$, or $E(N_{i_k})-E(N_{i_{k-1}})$ is a triad of $N_{i_k}$, when $k\ge 2$. 
 
The sequence of matroids that appear in the statement of the strong splitter theorem is
\begin{equation*}
N_0,N_{i_1},N_{i_2},\dots,N_{i_m},N_{i_m+1},N_{i_m+2},\dots,N_t
\end{equation*}
(That is, we remove all matroids having rank less than the rank of $M$ except for those having the rank for the first time.) $\qed$
\bigskip

\section {An application of the Strong Splitter Theorem}

In this section we will use the Strong Splitter Theorem to make some progress on determining the class of almost-regular matroids. A non-graphic matroid $M$ is   {\it almost-graphic} if, for all elements e, either $M\backslash e$ or $M/e$ is graphic. A non-regular matroid is   {\it almost-regular} if, for all                                             
elements $e$, either $M\backslash e$ or $M/e$ is regular. An element $e$ for which both                                             
$M\backslash e$ and $M/e$ are regular is called a {\it regular element}. Determining these classes of matroids was listed as an unsolved problem in the first edition of Oxley's book {\it Matroid Theory.}   In [4] we determined completely the class of almost-graphic matroids as well as the class of almost-regular matroids with at least one regular element. The problem that remains to be solved appears in the second edition as follows: [6, 15.9.8]: {\it Find all non-regular matroids $M$ such that, for all elements $e$, exactly one of $M\backslash e$ and $M/e$ is regular.}

In order to determine the class of almost-regular matroids with at least one regular element, we turned the problem into a series of excluded-minor classes and determined the members in them. In the almost-graphic paper, we proved the following characterization of almost-regular matroids with no $E_5$-minor [4, 8.2]. 
\bigskip
             
\noindent {\bf Theorem 3.1.} {\it Suppose $M$ is a $3$-connected binary almost-regular matroid with no $E_5$-minor.  Then $M\cong X_{12}$ or $M$ or $M^*$ is isomorphic to a $3$-connected restriction of $S_{3n+1}$ for $n\ge 3$, $\mathcal F_1(m, n, r)$ or $\mathcal F_2(m, n, r)$ for $m, n, r\ge 1$. }  
\bigskip

See [4] for a detailed description of the infinite families as well as the exceptional matroid $X_{12}$ that is a splitter for the class. Matroids like $T_{12}$ and $X_{12}$ and the rank-5, 10-element matroids $E_4$ and $E_5$ (which are single-element coextensions of $P_9$) play a useful role in the structure of binary matroids and feature in several papers [5]. Matrix representations for $E_4$ and $E_5$ are shown below.
\small 
\[
E_4=\left[ 
\begin{array}{c|ccccc}
&    0&1&1&1&1 \\
&    1&0&1&1&0 \\
I_5& 1&1&0&1&0 \\
&    1&1&1&1&0 \\
&    1&1&0&0&1
\end{array} 
\right] 
E_5=\left[ 
\begin{array}{c|ccccc}
&    0&1&1&1&1 \\
&    1&0&1&1&0 \\
I_5& 1&1&0&1&1 \\
&    1&1&1&1&0 \\
&    1&1&0&0&0
\end{array} 
\right] 
\]
\normalsize

To finish the almost-regular problem, one has to determine the almost-regular matroids with an $E_5$-minor. This is complicated by the presence of internally 4-connected members. However, in this paper we establish that the only thing left to do is to find the almost-regular matroids with both an $E_5$ and an $E_4$-minor.  

We found a rank-6, 12-element self-dual matroid that is a splitter for the class of binary 3-connected matroids with an $E_5$-minor and no $E_4$-minor. A matrix representations for $M_{12}$ is shown below.

\small
\[
M_{12}=\left[ 
\begin{array}{c|cccccc}
&     0 &  1  & 1 &  1 &  1 &  1  \\
&     1 &  0  & 1 &  1 &  0 &  0 \\
I_6&  1 &  1  & 0 &  1 &  1 &  0 \\
&     1 &  1  & 1 &  1 &  0 &  1  \\
&     1 &  1  & 0 &  0 &  0 &  1 \\
&     1 &  0  & 0 &  1 &  1 &  1 
\end{array}
\right] 
\] 
\normalsize

The main theorem in this section determines the class of matroids with an $E_5$-minor, but no $E_4$-minor. It has a finite list of members. With the exception of $M_{12}$, they are among the 3-connected restrictions of a rank-5, 17 element matroid, that we call $R_{17}$. The matroid $R_{17}$ is an extension of both $E_5$ and $R_{10}$ (the unique splitter for regular matroids). Note that $R_{17}$ has 3-connected restrictions that do not have an $E_5$-minor. A matrix representation for $R_{17}$ is shown below.  

\small
\[ 
R_{17}=\left[ 
\begin{array}{c|cccccccccccc}
&    0&1&1&1&1&0&0&0&0&1&1&1 \\
&    1&0&1&1&0&0&1&1&1&0&0&1 \\
I_5& 1&1&0&1&1&1&0&0&1&0&1&1 \\
&    1&1&1&1&0&1&0&1&0&1&0&0 \\
&    1&1&0&0&0&1&1&0&0&1&1&1
\end{array}
\right] 
\] 
\normalsize

\bigskip

\noindent The next result is the main theorem of this section.

\bigskip

\noindent {\bf Theorem 3.2.} {\it Suppose $M$ is a $3$-connected binary matroid with an $E_5$-minor and no $E_4$-minor.  Then $M\cong M_{12}$ or $M$ or $M^*$ is isomorphic to $R_{17}$ or is a $3$-connected restrictions of $R_{17}$ having an $E_5$-minor.}  
\bigskip

\noindent {\bf Proof.} The matroid $E_5$ is self-dual and has seven non-isomorphic binary 3-connected single-element extensions, shown in Appendix Table A1. Observe that all the extensions, except $A$, $B$ and $C$ have an $E_4$-minor. Matrix representations for $A$, $B$, and $C$ are given below.
\tiny
\[ 
A=\left[ 
\begin{array}{c|cccccc}
&    0&1&1&1&1&0 \\
&    1&0&1&1&0&0 \\
I_5& 1&1&0&1&1&1 \\
&    1&1&1&1&0&0 \\
&    1&1&0&0&0&1
\end{array} 
\right] 
B=\left[ 
\begin{array}{c|ccccccc}
&    0&1&1&1&1&1 \\
&    1&0&1&1&0&0 \\
I_5& 1&1&0&1&1&0 \\
&    1&1&1&1&0&1 \\
&    1&1&0&0&0&1
\end{array} 
\right] 
C=\left[ 
\begin{array}{c|ccccccc}
&    0&1&1&1&1&1 \\
&    1&0&1&1&0&1 \\
I_5& 1&1&0&1&1&0 \\
&    1&1&1&1&0&0 \\
&    1&1&0&0&0&1
\end{array} 
\right] 
\] 
\normalsize

The proof is in three stages. First, we will show that all the coextensions of $A$, $B$, and $C$ have an  $E_4$-minor with the exception of $M_{12}$. Suppose $M$ is a coextension of $A$, $B$, $C$. Then a partial matrix representation for $M$ is shown in Figure 1. 

\begin{figure}[h]
\centering
\epsfxsize 2in \epsfbox{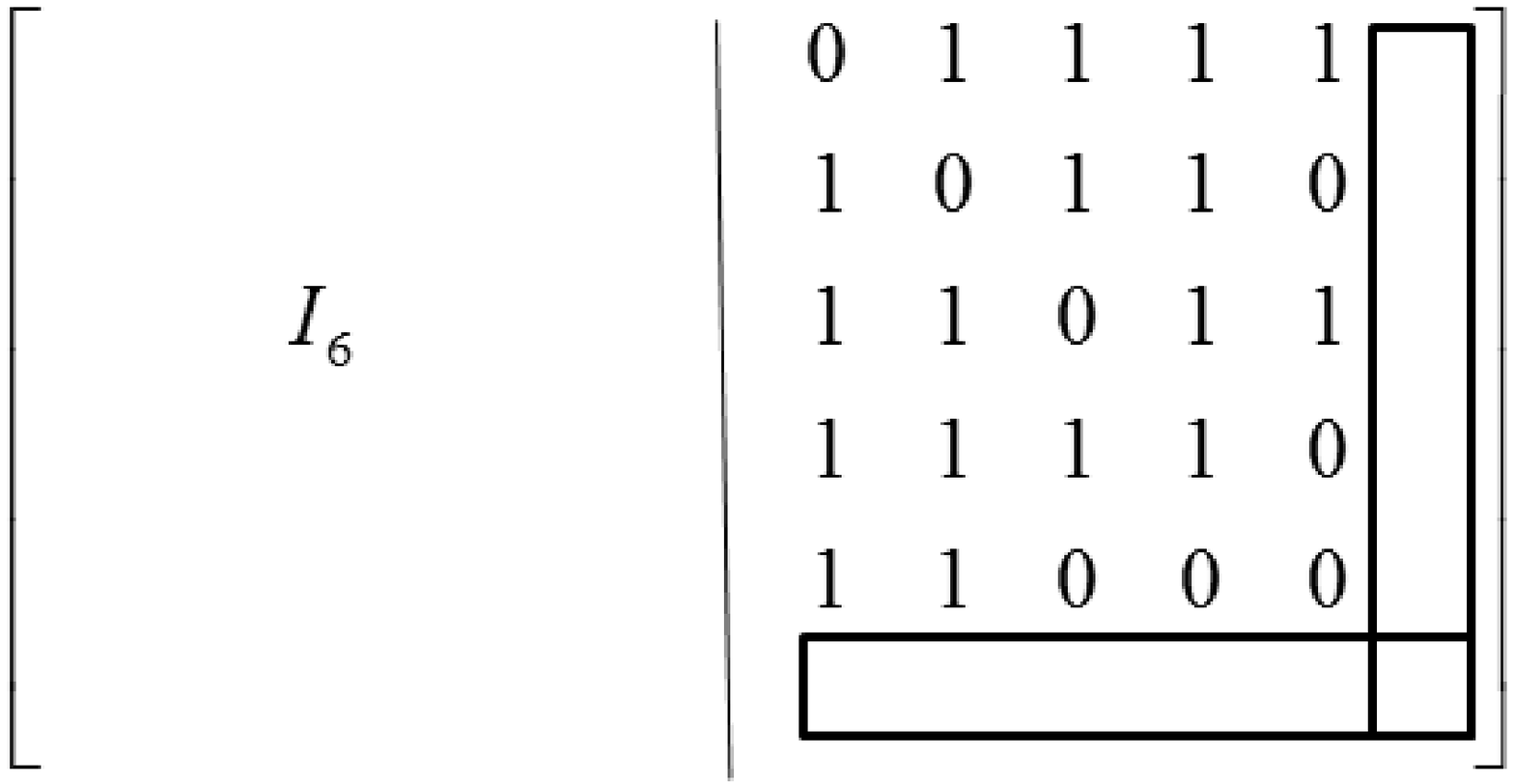}
\caption{Structure of a coextension of $A$, $B$, $C$ }
\end{figure}

There are three types of rows that may be inserted into the last row on the right-hand side of the matrix in Figure 1.

\begin {enumerate}
\item[(i)] rows that can be added to $E_5$ to obtain a coextension with no $E_4$-minor, with a 0 or 1 as the last entry; 
\item[(ii)] the identity rows with a 1 in the last position; 
 \item[(iii)] and the rows ``in-series" to the right-hand side of matrices $A$, $B$, $C$ with the last entry reversed. 
\end{enumerate}

\noindent Type I rows are $[0 0 1 1 1 0]$, $[0 0 1 1 1 1]$ $[0 1 0 0 1 0]$, $[0 1 0 0 1 1]$, $[0 1 0 1 0 0]$, $[0 1 0 1 0 1]$, $[0 1 1 0 0 0]$, $[0 1 1 0 0 1]$, $[1 0 0 1 1 0]$, $[1 0 0 1 1 1]$, $[1 0 1 0 1 0]$, $[1 0 1 0 1 1]$, $[1 1 1 0 1 0]$, and $[1 1 1 0 1 1]$. They are obtained from the Appendix Table A2.
Type II rows are $[1 0 0 0 0 1]$, $[0 1 0 0 0 1]$, $[0 0 1 0 0 1]$, $[0 0 0 1 0 1]$, and $[0 0 0 0 1 1]$. 
Type III rows are specific to the matrices $A$, $B$, $C$. For matrix $A$ they are $[0 1 1 1 1 1]$, $[1 0 1 1 0 1]$, $[1 1 0 1 1 0]$, $[1 1 1 1 0 1]$, $[1 1 0 0 0 0]$. For matrix $B$ they are $[0 1 1 1 1 0]$, $[1 0 1 1 0 1]$, $[1 1 0 1 1 1]$, $[1 1 1 1 0 0]$, and $[1 1 0 0 0 0]$. For $C$ they are $[0 1 1 1 1 0]$, $[1 0 1 1 0 0]$, $[1 1 0 1 1 1]$, $[1 1 1 1 0 1]$, and $1 1 0 0 0 0]$.

Most of the above rows result in matroids that have an $E_4$-minor. See red rows in the Appendix Table A3. Only a few coextensions must be specifically checked for an $E_4$-minor: $(A, coextn 11)$, $(B, coextn 8)$, $(C, coextn 8)$, $(C, coextn 9)$, $(C, coextn 10)$, $(C, coextn 12)$, and $(C, coextn 14)$. 
 
Observe that, 
$(A, coextn 11)/11\backslash 3 \cong E_4$, 
$(C, coextn 8)/12\backslash 2 \cong E_4$, 
$(C, coextn 9)/12\backslash 1 \cong E_4$, 
$(C, coextn 10)/12\backslash 10 \cong E_4$, and 
$(C, coextn 14)/12\backslash 6 \cong E_4$. 
Further, it is easy to check that $(B, coextn 8)\cong (C, coextn 12)$ and this matroid does not have an $E_4$-minor. This is the matroid $M_{12}$.

Second, we must establish that $M_{12}$ is a splitter for the class of matroids with an $E_5$-minor, but no $E_4$-minor. By Corollary 1.2 and the fact that $M_{12}$ is self-dual, we only need to check the single-element coextensions of $M_{12}$. From Appendix Table A3 observe that $M_{12}$, as a coextension of $C$, may be obtained by adding exactly one row. Thus there are no further rows that may be added to form coextensions without an $E_4$-minor. It follows that $M_{12}$ is a splitter for the class of binary matroids with an $E_5$, but no $E_4$-minor.

Third, we must show that if $M$ has an $E_5$ and no $E_4$-minor, then either $M\cong M_{12}$ or $r(M)\le 5$. To do this, let us begin by computing the single-element extensions of $A$, $B$, and $C$ with no $E_4$-minor. From Appendix Table A1, we may conclude that the only columns that can be added to $E_5$ to obtain a matroid with no $M^*(K_5\backslash e)$-minor are  $[0 0 1 0 1]$, $[0 0 1 1 0]$, $[0 1 0 1 1]$, $[0 1 1 0 0]$ $[1 0 0 1 1]$, $[1 1 0 0 1]$, $[1 1 1 0 1]$. Adding these columns gives us four non-isomorphic single-element extensions of $A$, $B$, and $C$. They are $D$, $E$, $F$, and $G$ shown below. 

\tiny
\[ 
D=\left[ 
\begin{array}{c|ccccccc}
&    0&1&1&1&1&0&0 \\
&    1&0&1&1&0&0&0 \\
I_5& 1&1&0&1&1&1&1 \\
&    1&1&1&1&0&0&1 \\
&    1&1&0&0&0&1&0
\end{array} 
\right] 
E=\left[ 
\begin{array}{c|ccccccc}
&    0&1&1&1&1&0&0 \\
&    1&0&1&1&0&0&1 \\
I_5& 1&1&0&1&1&1&0 \\
&    1&1&1&1&0&0&1 \\
&    1&1&0&0&0&1&1
\end{array} 
\right] 
\] 
\normalsize

\tiny
\[ 
F=\left[ 
\begin{array}{c|ccccccc}
&    0&1&1&1&1&0&1 \\
&    1&0&1&1&0&0&1 \\
I_5& 1&1&0&1&1&1&0 \\
&    1&1&1&1&0&0&0 \\
&    1&1&0&0&0&1&1
\end{array} 
\right] 
G=\left[ 
\begin{array}{c|ccccccc}
&    0&1&1&1&1&0&1 \\
&    1&0&1&1&0&0&1 \\
I_5& 1&1&0&1&1&1&1 \\
&    1&1&1&1&0&0&0 \\
&    1&1&0&0&0&1&1
\end{array} 
\right] 
\] 
\normalsize

\noindent Suppose $M$ is a coextension of $D$, $E$, $F$, or $G$. Then the structure of $M$ is shown in Figure 2. Observe that one row and two columns may be added to $E_5$.

\begin{figure}[h]
\centering
\epsfxsize 2in \epsfbox{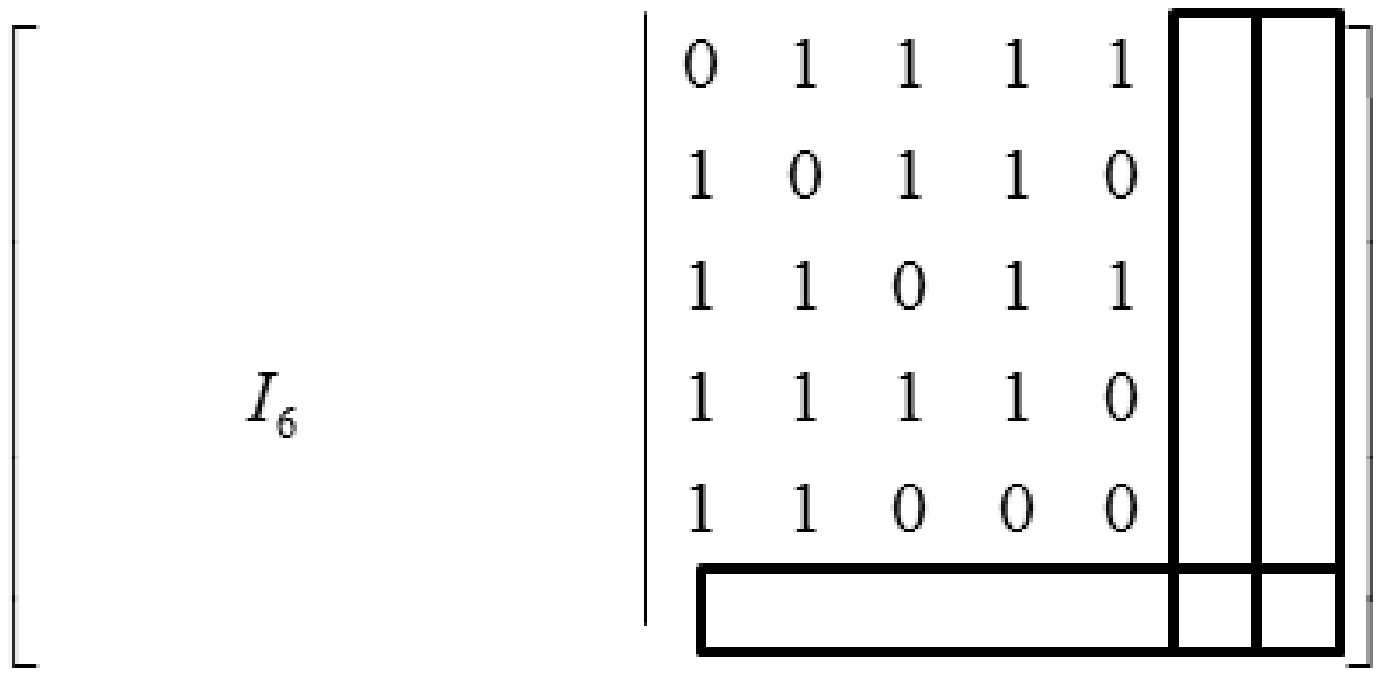}
\caption{Structure of a coextension of $D$, $E$, $F$, $G$ }
\end{figure}

\noindent Once again three types of rows may be added.

\begin {enumerate}
\item[(i)] the  rows that can be added to $D$, $E$, $F$, or $G$ to obtain a coextension with no $E_4$-minor, with a 0 or 1 in the last entry.   
\item[(ii)] the identity rows with a 1 in the last position; 
 \item[(iii)] and the rows ``in-series" to the right-hand side of the matrices with the last entry reversed. 
\end{enumerate}

\noindent Suppose $M$ is the coextension obtained by adding a Type I row, then $M\backslash 13$ is 3-connected. However, the only rank-6, 12-element matroid in the class is $M_{12}$ and it is a splitter; a contradiction. Thus we may assume $M\backslash 13$ is not 3-connected. 

 Adding a Type II or III row (with the exception of $[0 0 0 0 0 1 1]$) causes $M\backslash e$ to be 3-connected where $e\in \{12, 13\}$ (and again there are no such matroids except $M_{12}$ which is a splitter).  So the  only coextension we must check is the one formed by adding row $[0 0 0 0 0 1 1]$. That is the coextension in which $\{6, 11, 12\}$ is a triad. Let $D'$, $E'$, $F'$, and $G'$ be the coextensions of $D$, $E$, $F$, and $G$, respectively, obtained by coextending by row $[0 0 0 0 0 1 1]$. Then in each case we can find an $E_4$ minor. In particular, 
$D'/1\backslash \{3, 11\} \cong E_4$,  
$E'/1\backslash \{7, 11\} \cong E_4$, 
$F'/1\backslash \{7, 11\} \cong E_4$, and
$G'/1\backslash \{7, 11\} \cong E_4$.  

Finally, observe that if $M$ is an extension of $E_5$ of size $k\ge 13$, then for some $e\in \{11, \dots , k\}$, $M\backslash e$ is 3-connected. The result follows from Theorem 1.4. $\qed$

\bigskip

The next theorem uses Theorem 3.2 to characterize the almost-regular matroids with an $E_5$-minor, but no $E_4$-minor.    
\bigskip

\noindent {\bf Theorem 3.3.} {\it Suppose $M$ is a $3$-connected binary almost-regular matroid with an $E_5$-minor.  Then, either $M$ has an $E_4$-minor or $M\cong E_5$, $B$, or $B^*$. }  
\bigskip

\noindent {\bf Proof.} We begin by looking at the single-element extensions of $E_5$ and determining that $B$ and $H$ are the only ones that are almost-regular. Note that $B$ and $H$ are $E_{5, 11}$ and $D_{5, 11}$ in Appendex Table III in [4]. Further, note that $H$ has an $E_4$-minor. Since $B$ is formed by adding just one column to $E_5$, no further extension of $E_5$ is almost-regular without an $E_4$-minor. The result then follows from Theorem 3.2 because $M_{12}$ is not almost-regular. $\qed$
\bigskip

Thus we may assume an almost-regular matroid with an $E_5$-minor must also have an $E_4$-minor (with the exception of $B$ and $B^*$) and more importantly, it must have an $H$-minor, where $H$ is the matroid shown below:

\small
\[
H=\left[ 
\begin{array}{c|cccccc}
&     0 &  1  & 1 &  1 &  1 &  1  \\
&     1 &  0  & 1 &  1 &  0 &  0 \\
I_5&  1 &  1  & 0 &  1 &  1 &  1 \\
&     1 &  1  & 1 &  1 &  0 &  0  \\
&     1 &  1  & 0 &  0 &  0 &  1 
\end{array}
\right] 
\]
\normalsize

\noindent As a consequence of the above results we can strengthen  Theorem 3.1 as follows:   
\bigskip

\noindent {\bf Theorem 3.5.} {\it Suppose $M$ is a $3$-connected binary almost-regular matroid with no $H$ or $H^*$-minor.  Then $M\cong E_5$, $B$, $B^*$, $X_{12}$, or $M$ or $M^*$ is isomorphic to a $3$-connected restriction of $S_{3n+1}$ for $n\ge 3$, $\mathcal F_1(m, n, r)$ or $\mathcal F_2(m, n, r)$ for $m, n, r\ge 1$. }  
\bigskip

\bigskip

\bigskip
\noindent {\bf References}

\bigskip

\begin{enumerate}

\item R. E. Bixby and C. R. Coullard (1987).  Finding a 3-connected minor maintaining a fixed minor and a fixed element. {\it Combinatorica} {\bf 7}, 231-242.

\item T. H. Brylawski (1972). A decomposition for combinatorial geometries. {\it Transactions of the Americal Mathematical Society} {\bf 171}, 235-282.

\item C. R. Coullard and J. G. Oxley (1992). Extension of Tutte's wheels and whirls theorem. {\it Journal of Combinatorial Theory Series B}  {\bf 56}, 130 - 140.

\item S. R. Kingan  and M. Lemos  (2002), Almost-graphic matroids. {\it Advances in Applied Mathematics} {\bf 28}, 438 - 477.

\item S. R. Kingan  and M. Lemos  (to appear), Matroids with at least two regular elements. {\it European Journal of Combinatorics}.

\item J. G. Oxley (1992). {\it Matroid Theory}. Oxford University Press, New York.

\item P. D. Seymour (1980). Decomposition of regular matroids.  {\it Journal of Combinatorial Theory Series B } {\bf 28 }, 305-359.   

\item J. J. -M. Tan (1981). Matroid 3-connectivity. Ph.D. thesis, Carleton University.

\item K. Truemper (1992). {\it Matroid Decomposition}. Academic Press, Boston.   

\item W. T. Tutte (1960). An algorithm for determining whether a given binary matroid is graphic. {\it Proceedings of the American Mathematical Society} {\bf 11}, 905-917.

\item W. T. Tutte (1966). Connectivity in matroids. {\it Canadian Journal of Mathematics} {\bf 18}, 1301-1324.
\end {enumerate}

\bigskip

\noindent  { \bf Appendix}  
\bigskip

\tiny
 \begin{center}
\begin{tabular}{|p{15em}|c|c|}
\hline
 \bf{Extension Columns} & {\bf Name} & {\bf $E_4$-minor}  \\  \hline \hline
  $[0 0 1 0 1]$ $[0 0 1 1 0]$ $[0 1 0 1 1]$ $[0 1 1 0 0]$  & $A$ &  No   \\  \hline

 $[1 0 0 1 1]$ & $B$   & No   \\  \hline

 $[1 1 0 0 1]$ $[1 1 1 0 1]$   &   $C$   & No   \\  \hline

 $[0 0 0 1 1 ]$ $[0 0 1 1 1]$ $[0 1 0 0 1]$ $[0 1 1 0 1]$  && Yes   \\  \hline

 $0 1 0 1 0]$ $[0 1 1 1 0]$ && Yes     \\  \hline

 $[1 0 0 0 1]$ $[1 0 0 1 0]$  $[1 1 0 1 1]$ $[1 1 1 0 0]$  && Yes   \\  \hline

 $[1 0 1 0 1]$ $[1 0 1 1 0]$ $[1 1 0 0 0]$ $[1 1 1 1 1]$   &H& Yes   \\  \hline
\end{tabular}
 \end{center}
\normalsize
 \begin{center} \small Table A1: Single-element extensions of $E_5$ \normalsize \end{center} 

\bigskip

\tiny
 \begin{center}
\begin{tabular}{|p{15em}|c|c|c|}
\hline
\bf{Coextension Rows} & {\bf Name} \\      \hline \hline
$[0 0 1 1 1]$ $[0 1 0 0 1]$ $[0 1 0 1 0]$ $[0 1 1 0 0]$  & $A^*$  \\  \hline

$[1 0 0 1 1]$ & $B^*$      \\  \hline

$[1 0 1 0 1]$ $[1 1 1 0 1]$   &   $C^*$    \\  \hline

$[0 0 0 1 1 ]$ $[0 0 1 0 1]$ $[0 1 0 1 1]$ $[0 1 1 0 1]$   & \\  \hline

$[0 0 1 1 0]$ $[0 1 1 1 0]$   &   \\  \hline

$[1 0 0 0 1]$ $[1 0 0 1 0]$  $[1 0 1 1 1]$ $[1 1 1 0 0]$ & \\  \hline

$[1 0 1 0 0]$ $[1 1 0 0 1]$ $[1 1 0 1 0]$ $[1 1 1 1 1]$ & $H^*$ \\  \hline
\end{tabular}
 \end{center}
\normalsize
 \begin{center} \small Table A2: Single-element coextensions of $E_5$ \end{center} 
\normalsize

\tiny
 \begin{center}
\begin{tabular}{|c|c|p{35em}|}
\hline
\bf{Matroid} & Name  & {\bf Coextension Row}  \\  \hline \hline
$A$ & coext 1 & \color {red} $[0 0 0 0 1 1]$  $[0 0 0 1 0 1]$ \color {black} $[0 0 1 0 1 0]$ $[0 1 1 0 1 0]$ $[1 0 1 1 1 1]$ $[1 1 1 0 0 1]$               \\  \hline

& coext 2 & $[0 0 0 1 1 0]$  $[1 1 0 0 1 1]$  $[1 1 0 1 0 1]$    \\  \hline

& coext 3 & $[0 0 0 1 1 1]$ \color {red} $[1 0 1 0 1 1]$ $[1 1 1 0 1 1]$ \\  \hline

& coext 4 & \color {red}$[0 0 1 0 0 1]$\color {black} $[0 1 0 1 1 0]$ \color {red}$[0 1 1 1 1 1]$ \color {black}      \\ \hline

& coext 5 & $[0 0 1 0 1 1]$ $[0 1 1 0 1 1]$ \color {red} $[1 0 0 1 1 1]$      \\ \hline

& coext 6 & $[0 0 1 1 0 0]$ $[0 1 1 1 0 0]$ \color {red}$[1 1 0 0 0 0]$      \color {black} \\ \hline

& coext 7 & $[0 0 1 1 0 1]$ \color {red} $[0 1 0 0 1 0]$ $[0 1 0 1 0 0]$  \color {black}  $[0 1 1 1 0 1]$  $[1 0 1 1 1 0]$  $[1 1 1 0 0 0]$     \\ \hline

& coext 8 & \color {red} $[0 0 1 1 1 0]$ $[0 1 1 0 0 0]$  $[1 0 1 1 0 1]$ \color {black} $[1 1 0 0 1 0]$  $[1 1 0 1 0 0]$ \color {red} $[1 1 1 1 0 1]$      \\ \hline

& coext 9 & \color {red}  $[0 0 1 1 1 1]$ $[0 1 1 0 0 1]$ \color {black} $[1 0 0 0 1 1]$  $[1 0 0 1 0 1]$ \color {red}  $[1 0 1 0 1 0]$  $[1 1 1 0 1 0]$     \\ \hline

& coext 10 & \color {red} $[0 1 0 0 0 1]$   $[1 0 0 0 1 0]$ \color {black} $[1 0 0 1 0 0]$     \\ \hline

& coext 11 & \color {red}  $[0 1 0 0 1 1]$ $[0 1 0 1 0 1]$ $[1 0 0 1 1 0]$    \\ \hline

& coext 12 & $[0 1 0 1 1 1]$ \\ \hline

& coext 13 & \color {red} $[1 0 0 0 0 1]$ \color {black}  $[1 0 1 0 0 0]$ $[1 1 1 1 1 0]$    \\ \hline

& coext 14 & $[1 0 1 0 0 1]$ \color {red} $[1 1 0 1 1 0]$ \color {black} $[1 1 1 1 1 1]$     \\ \hline

 \hline

$B$ & coext 1  & \color {red}$[0 0 0 0 1 1]$ $[0 0 0 1 0 1]$ \color {black}  $[0 0 0 1 1 0]$  \color {red} $[0 0 1 0 0 1]$ \color {black}  $[0 0 1 0 1 0]$ \color {red} $[0 0 1 1 1 1]$ $[0 1 0 0 1 0]$ $[0 1 0 1 0 0]$ \color {black}  $[0 1 0 1 1 1 ]$ \color {red}  $[0 1 1 0 0 0]$ \color {black}  $[0 1 1 0 1 1]$ \color {red} $[0 1 1 1 1 0]$ \color {black}             \\  \hline

& coext 2 &  $[0 0 0 1 1 1]$ $[0 0 1 0 1 1]$ $[0 1 0 1 1 0]$ $[0 1 1 0 1 0]$            \\  \hline

& coext 3 &  $[0 0 1 1 0 0]$ \color {red} $[0 1 0 0 0 1]$ \color {black}  $[0 1 1 1 0 1]$           \\  \hline

& coext 4 &  $[0 0 1 1 0 1]$ \color {red}  $[0 0 1 1 1 0]$ $[0 1 0 0 1 1]$ $[0 1 0 1 0 1]$ $[0 1 1 0 0 1]$ \color {black}  $[0 1 1 1 0 0]$             \\   \hline

& coext 5 &  \color {red} $[1 0 0 0 0 1]$ \color {black} $[1 0 0 0 1 0]$ $[1 0 0 1 0 0]$ $[1 0 1 0 0 0]$ \color {red} $[1 0 1 1 0 1]$ \color {black} $[1 0 1 1 1 0]$ \color {red} $[1 1 0 0 0 0]$ \color {black} $[1 1 0 0 1 1]$ $[1 1 0 1 0 1]$ $[1 1 1 0 0 1]$ \color {red} $[1 1 1 1 0 0]$ \color {black} $[1 1 1 1 1 1]$             \\   \hline

& coext 6 & $[1 0 0 0 1 1]$ $[1 0 0 1 0 1]$ \color {red}  $[1 0 1 0 1 0]$ \color {black}  $[1 0 1 1 1 1]$ $[1 1 1 0 0 0]$ \color {red} $[1 1 1 0 1 1]$              \\   \hline

& coext 7 & \color {red}  $[1 0 0 1 1 0]$ \color {black} $[1 0 1 0 0 1]$ $[1 1 0 0 1 0]$ $[1 1 0 1 0 0]$ \color {red} $[1 1 0 1 1 1]$ \color {black} $[1 1 1 1 1 0]$            \\   \hline

& coext 8 & \color {red}  $[1 0 0 1 1 1]$ $[1 0 1 0 1 1]$ $[1 1 1 0 1 0]$             \\   \hline\hline

$C$ & coext 1 & \color {red} $[0 0 0 0 1 1]$ $[0 0 0 1 0 1]$ $[0 0 1 0 0 1]$  $[0 0 1 1 1 1]$         $[0 1 0 0 1 0]$ $[0 1 0 1 0 0]$  $[0 1 1 0 0 0]$  $[0 1 1 1 1 0]$ \color {black} $[1 0 0 0 1 0]$ $[1 0 0 1 0 0]$ $[1 0 1 0 0 0]$ $[1 0 1 1 1 0]$ $[1 1 0 0 1 1]$ $[1 1 0 1 0 1]$ $[1 1 1 0 0 1]$ $[1 1 1 1 1 1]$                   \\  \hline

& coext 2 & $[0 0 0 1 1 0]$  $[0 1 0 1 1 1]$     \\  \hline

& coext 3 & $[0 0 0 1 1 1]$  $[0 1 0 1 1 0]$  \color {red} $[1 0 0 1 1 0]$  $[1 1 0 1 1 1]$ \color {black}  \\  \hline

& coext 4 & $[0 0 1 0 1 0]$  $[0 1 1 0 1 1]$    \\  \hline

& coext 5 & $[0 0 1 0 1 1]$  $[0 1 1 0 1 0]$ \color {red} $[1 0 1 0 1 0]$ $[1 1 1 0 1 1]$     \\  \hline

& coext 6 & $[0 0 1 1 0 0]$  $[0 1 1 1 0 1]$     \\  \hline

& coext 7 & $[0 0 1 1 0 1]$  $[0 1 1 1 0 0]$ \color {red} $[1 0 1 1 0 0]$ $[1 1 1 1 0 1]$ \color {black}  \\  \hline

& coext 8 & \color {red} $[0 0 1 1 1 0]$  $[0 1 0 0 1 1]$ $[0 1 0 1 0 1]$ $[0 1 1 0 0 1]$     \\  \hline

& coext 9 & \color {red} $[0 1 0 0 0 1]$ \color {black}   \\  \hline

& coext 10 & \color {red} $[1 0 0 0 0 1]$  $[1 1 0 0 0 0]$ \color {black}      \\  \hline

& coext 11 & $[1 0 0 0 1 1]$  $[1 0 0 1 0 1]$ $[1 0 1 1 1 1]$ $[1 1 1 0 0 0]$    \\  \hline

& coext 12 & \color {red}  $[1 0 0 1 1 1]$     \\  \hline

& coext 13 & $[1 0 1 0 0 1]$  $[1 1 0 0 1 0]$ $[1 1 0 1 0 0]$ $[1 1 1 1 1 0]$   \\  \hline

& coext 14 & \color {red}  $[1 0 1 0 1 1]$  $[1 1 1 0 1 0]$  \\  \hline

\hline
\end{tabular}
\end{center}
\normalsize
 \begin{center} Table A3: Single-element coextensions of $A$ $B$ and $C$ \end{center}

\end {document}